\documentclass{amsart}
\usepackage{amsmath}
  \usepackage{paralist}
  \usepackage{graphics} 
  \usepackage{epsfig} 
 \usepackage[colorlinks=true]{hyperref}
\hypersetup{urlcolor=blue, citecolor=red}

\theoremstyle{definition}

\title[Prolonging Delay in CTL Activation]
      {Do Longer Delays Matter? \\ The Effect of Prolonging Delay in CTL Activation}

\author[Angela  Gallegos and Ami Radunskaya]{}

\subjclass{Primary: 34K60, 37N25; Secondary: 92B05.}
 \keywords{Delay Differential Equations, Stability Analysis.}

 \email{angela@oxy.edu}
 \email{aer04747@pomona.edu}

\begin{document}
\maketitle

\centerline{\scshape Angela Gallegos }
\medskip
{\footnotesize
 \centerline{Department of Mathematics, Occidental College}
   \centerline{ Los Angeles, CA USA}
} 

\medskip

\centerline{\scshape Ami Radunskaya}
\medskip
{\footnotesize
 \centerline{ Department of Mathematics, Pomona College}
   \centerline{Claremont, CA  USA}
}

\bigskip


\begin{abstract}
The activation of a specific immune response takes place in the lymphoid organs such as the spleen.  We present here a simplified model of the proliferation of specific immune cells in the form of a single delay equation.  We show that the system can undergo switches in stability as the delay is increased, and we interpret these results in the context of sustaining an effective immune response to a dendritic cell vaccine.
\end{abstract}

\section{Introduction}

Immunotherapy seeks to arm the patient with an increased number of effective immune cells that are specifically trained to seek out and destroy harmful cells.  One target for immunotherapy that has earned a place in the clinical spotlight recently is melanoma, since specific peptides have been isolated from melanoma cells that can be used in cancer vaccines.  A subset of immune cells, called dendritic cells, are able to stimulate the production of immune effector cells that are capable of recognizing and killing specific tumor cells.  Researchers have been able to extract these cells from patients and culture them {\it ex vivo} to create a vaccine that can boost a patient's response against their own cancerous cells, \cite{PilonThomas2004}.  The success of clinical trials of dendritic cell (DC) vaccines has resulted in the recent FDA approval of the first cancer vaccine for prostate cancer, \cite{ABCNews}

Despite promising clinical responses in vaccine trials, it remains difficult to predict which patients will actually respond to these vaccines and why,
\cite{Boon2006, {Trefzer2005}}.  This may be due to the complicated kinetics of the immune response to the presentation of antigen by the dendritic cells.  In this paper, we analyze the effect of the delay between the time the DCs come into contact with the T-cell population and the initiation of the expansion, or proliferating phase.  We focus on the dynamics describing the active effector T-cell population in the spleen, using a phenomenological description that captures the effect of the delay.  Our goal here is to isolate the role of this delay term in order to analyze its effect on the sustainability of T-cell production.

The model we present is very simple, and is not intended to realistically capture the entire cascade of immune events in the response to antigen presentation by dendritic cells.  Nor does this model describe the trafficking of immune cells between compartments of the body, in particular between lymphoid organs, such as the spleen, and the tumor site.  These effects are crucial to the understanding of the immune response and must be considered in the design of vaccine therapies. One group of experiments that looks at the trafficking of DCs between compartments is described in \cite{Preynat2007}, and a mathematical model of DC trafficking  and T-cell activation is presented in \cite{Ludewig2004}.    Elsewhere we combine these ideas, along with previous work on cell-lysis rates, \cite{dePillis2005}, in a more realistic model that includes several types of immune cells, and trafficking between compartments.

In section \ref{ModelDescription} we describe the model itself, and then we show stability switching as a function of the delay in Section \ref{Stability}.  We conclude with some discussion of the model simulations and results in Section \ref{Conclusion}.

\section{Description of the Mathematical Model} \label{ModelDescription}

\begin{figure}
\includegraphics[scale = .33]{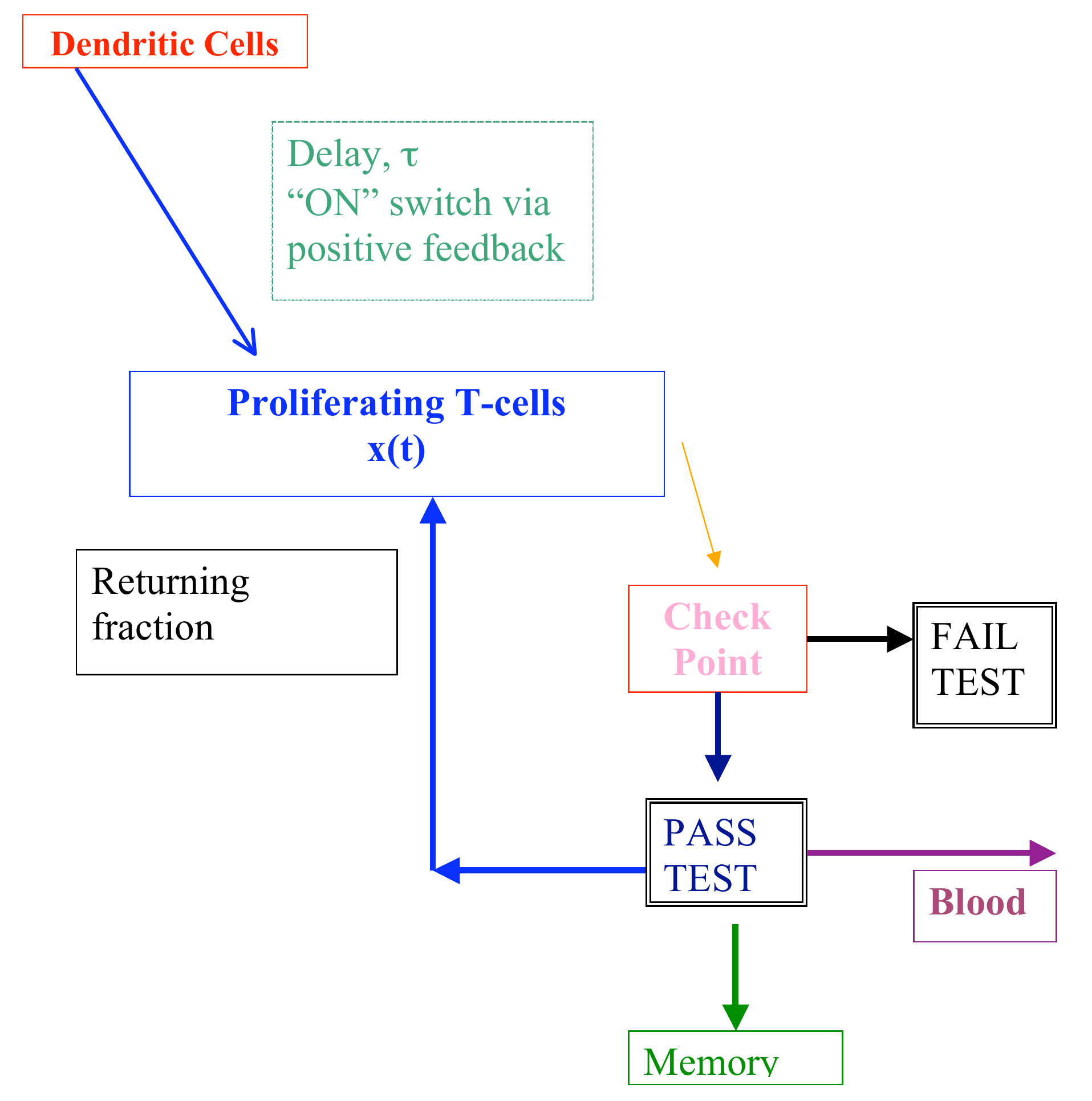}
\caption{A graphical representation of the process underlying the mathematical model.} \label{Fig:ModelProcess}
\end{figure}

After encountering antigen, such as tumor cells or tumor peptides, dendritic cells migrate to lymph organs such as the spleen.  As they migrate, they mature and, upon arrival in the spleen, they are able to activate the proliferation of T-cells that have the ability to seek out and destroy (\emph{lyse}) the specific target cells that provided the antigen.  This activation requires some contact time between the dendritic cells and the naive T-cells; this connection time is called the ``synaptic connection time." In our model, we represent this time as a delay, $\tau$.  Once proliferation is initiated, new antigen-specific T-cells are created.  A somewhat simplified description of the fate of these new cells is as follows:  The cells must go through a check-point during which the new cells are tested for efficacy against the antigen.  If they pass this test, they either move into the blood stream to migrate to the target cells, they become memory cells in order to protect against future challenges by the same antigen, or they return to the proliferating compartment to produce more specific T-cells.  Proliferation continues only while certain cytokines are present.  These cytokines are produced by previously activated T-cells that have been in contact with dendritic cells for the synaptic connection time, $\tau$.  When the dendritic cells have done their job ``educating" the T-cells, it is believed that they die, \cite{PilonThomas2004}, and are then ``cleared" from the spleen compartment.  When the DCs are no longer present, cell proliferation ceases, and the activated T-cells either move into the blood stream or become memory cells.  This process is represented graphically in Figure \ref{Fig:ModelProcess}.

In this model, the concentration of the dendritic cells at a given time, $t$, is given by $a(t)$, while the number of activated T-cells is given by $x(t)$.  The activation/proliferation rate is then given by $$a(t - \tau) f(x(t - \tau) ) x(t),$$ where $f$ is represents the ``feedback" function.  The ``feedback" function depends on the number of T-cells  present in the proliferating compartment  time $\tau$ earlier in the presence of the DCs; therefore they are producing the necessary cytokines. The current number of proliferating T-cells is represented by $x(t)$.   We assume that the ``feedback" function, $f$, is a positive function that increases to some maximum level, $r$, and then decreases.  The justification for this functional form is that activation by DCs is maximized at a certain concentration of DCs in the spleen.  The addition of more DCs has the effect of keeping the activated cells in the area.  If they are present for too long, the effectiveness of the T-cells begins to diminish.  This is sometimes called the ``decreasing potential hypothesis" in the literature, \cite{Joshi2008, Kalia2006}.  One simple function that captures these dynamics is given by:

\begin{center}
 \vskip .3in \hbox{\hskip .5in$ f(x_{\tau} ) = \displaystyle{\frac{  r x_{\tau}}{1 + x_{\tau}^4} } $ } \vskip -.5in \hskip .2in  \hbox{\includegraphics[scale = .15]{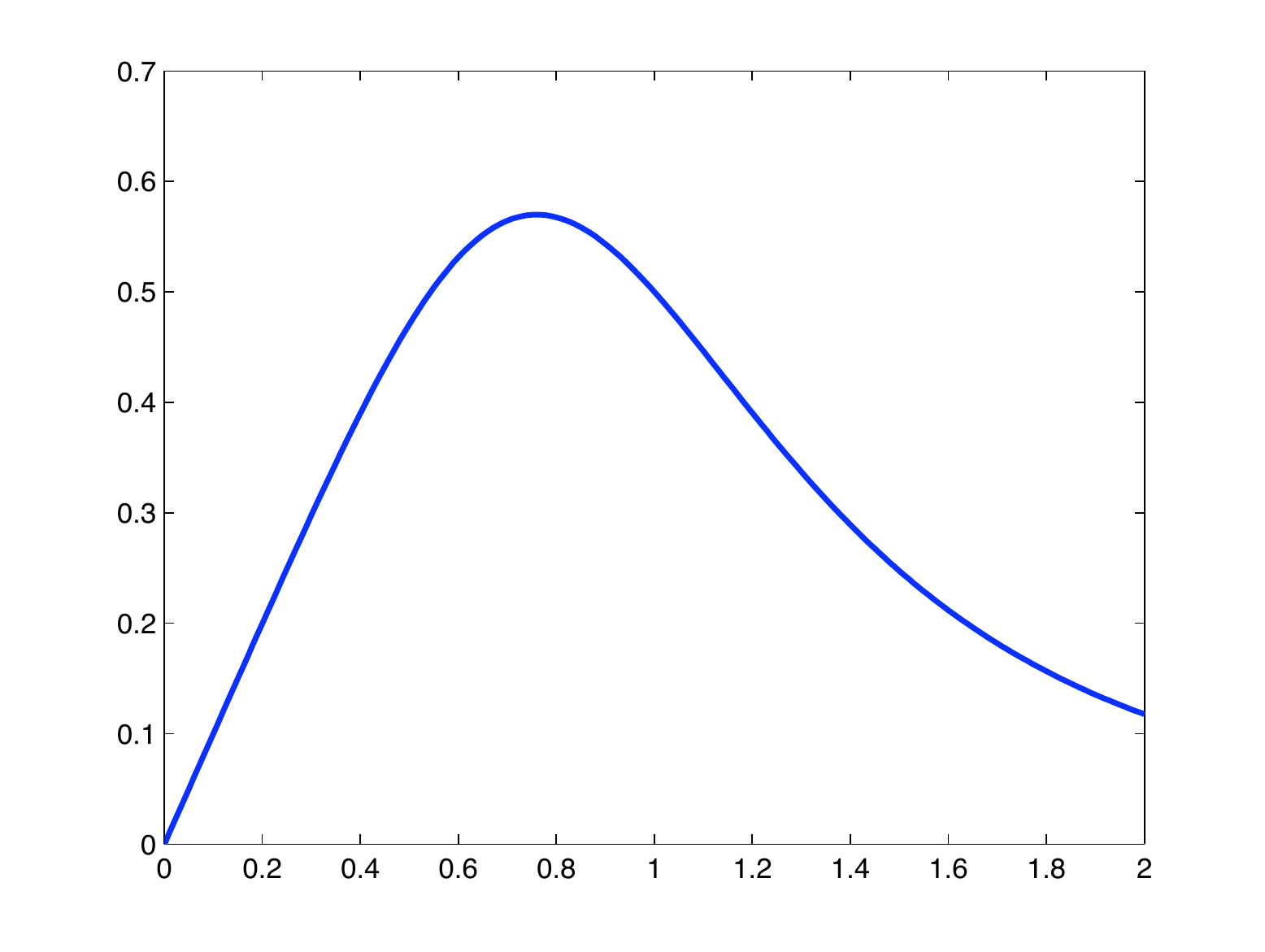} }
\end{center}
where $x_{\tau} = x(t - \tau)$ represents the population of active T-cells from time $\tau$ ago. 

We assume that T-cells leave the proliferating compartment  through death, moving into the blood stream or the memory compartment at an average per-cell rate, $\mu > 0 $.   Putting these pieces together results in the following single delay equation:
\begin{equation} \label{Eqn:DDE}
 \frac{dx}{dt} = a(t - \tau ) \,  f(x_{\tau}) x   - \mu x 
 \end{equation} 
where for simplicity:
$$ a(t) = \begin{cases} 1 & \hbox{ when DCs are present} \\ 0 & \hbox{otherwise} \end{cases} .$$

By rescaling the time variable, we can assume that $r = 1$, (the decay parameter, $\mu$, must also be appropriately scaled). Simulated solutions are shown in Figures \ref{Fig:SampleTrajectories_tauzero} and \ref{Fig:SwitchingTimes}.  These solutions show that changing the value of $\tau$, as well as the initial history (the values of $x(t)$ for $- \tau < t < 0$), can have a strong effect on the long-term behavior of the system.  We verify this in the next section by performing a stability analysis.

\section{Stability Anlaysis} \label{Stability}

When dendritic cells are present in the spleen, $a(t) = 1$, so that the differential equation becomes:
\begin{equation} \label{Eqn:ScaledDDE}
 \frac{dx}{dt} = \frac{x_{\tau} }{1 + x_{\tau}^4 } x  - \mu x.
 \end{equation} 

\subsection{Equilibria and stability when $\tau = 0$.}

\begin{figure}
\includegraphics[scale = .5]{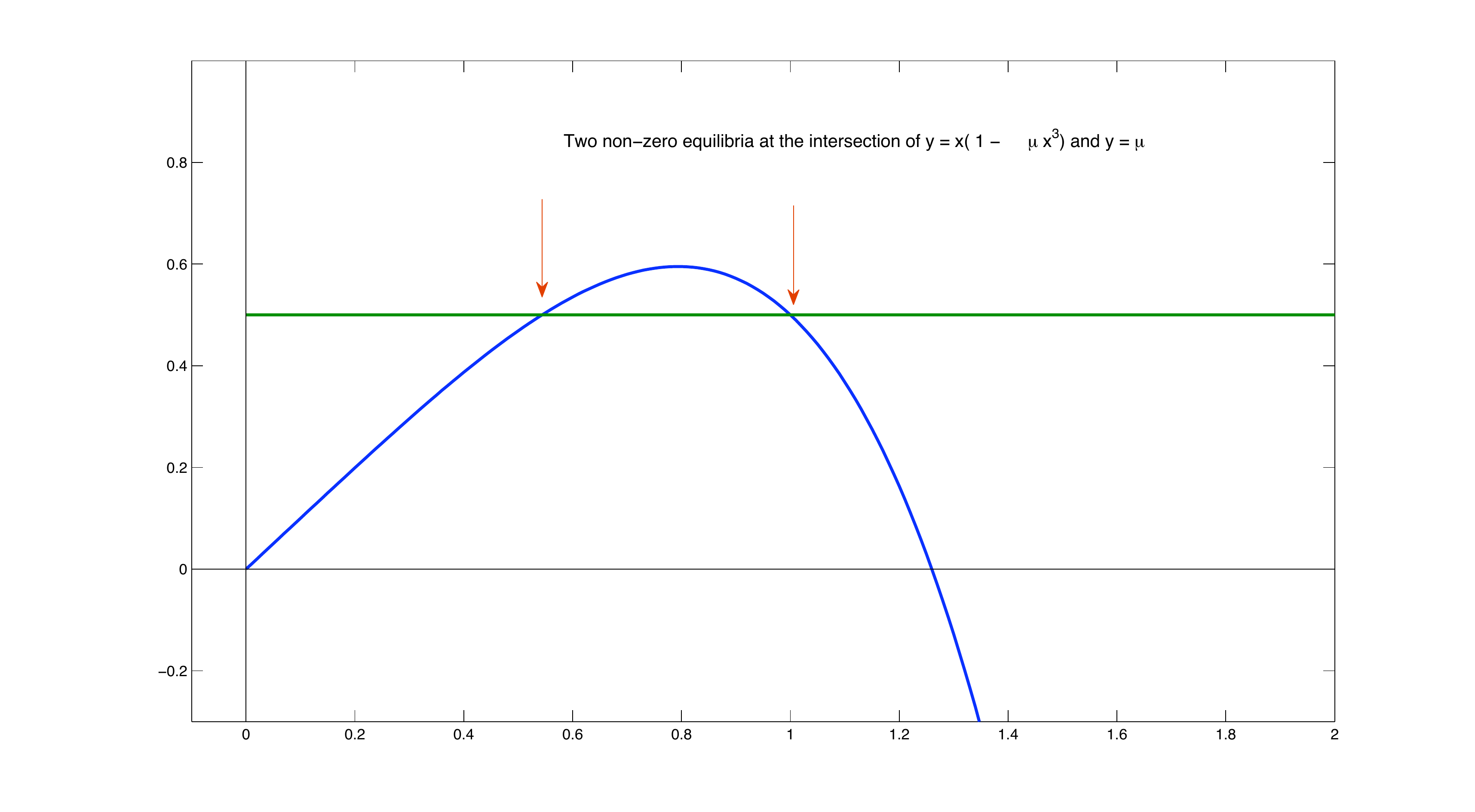}
\caption{The non-zero equilibria of the DDE given in \ref{Eqn:ScaledDDE} are the solutions of
$x (1 - \mu x^3) = \mu$.  These are positive for $\mu$-values in $(0, \frac{3^{3/4}}{4})$.}
\label{Fig:TwoEquil}
\end{figure}

Equation \ref{Eqn:ScaledDDE} potentially has three equilibria.  One is the zero fixed point: $x_E = 0$, and the other two are solutions to
$$ \frac{x}{1 + x^4} = \mu \Leftrightarrow x ( 1 - \mu x^3) = \mu .$$
These two solutions are depicted graphically in Figure \ref{Fig:TwoEquil}. Note that the function $x   - \mu x^4$ has a maximum which occurs at $x_{Max} = (4\mu)^{-1/3}$.  When the maximum of this function is greater than $\mu$ the two additional equilibria will exist.  This occurs for $\mu \in \left({ 0, \displaystyle{\frac{3^{3/4}}{4}} }\right)$.    We denote the greater of these two equilibria by $x^*$, and the smaller by $x^-$, if they exist.

If $x^*$ is stable, then proliferation is sustained, and active effector cells continue to be produced and sent into the blood stream and into the memory compartment.  However, if $x^*$ is unstable,  then production may die off and no new effector cells will be produced, even when DCs are still present in the spleen.
Our goal in this section is to show that changing the value of the delay, $\tau$, can change the stability of $x^*$.  

Linearization around an equilibrium, $x_E$ gives:
$$ \frac{dz}{dt} = \alpha  z_\tau + \beta z, $$
where $z = x - x_E $, \quad $\alpha = x \frac{\partial f}{\partial x_{\tau}} $, evaluated at the equilibrium $x_E$, and  $ \beta =  f(x_E) - \mu$.  The characteristic equation is:
$$ \lambda = \alpha  e^{- \lambda \tau} + \beta.$$
Evaluating $\alpha$ and $\beta$ at $x = x_E$ gives:
$$ \alpha = \frac{ \mu^2}{x_E}  \left({ 1 - 3 x_E^4}\right), \quad \hbox{ when } x_E = x^*, x^-, \quad \hbox{and}  \quad \alpha = 0 \hbox{ when } \, x_E = 0 $$
$$ \beta = 0, \quad x_E = x^*, x^-, \quad \hbox{and} \quad \beta = - \mu, \, x_E = 0 .$$
When we consider the characteristic equation in the case of the equilibrium $x_E = 0$, we always obtain $\lambda = -\mu < 0$, \emph{independent of $\tau$.}  Thus we see that the equilibrium at $x_E =0$ is {\bf always stable}.

On the other hand,linearizing around the upper equilibrium at $x^*$ in the case when $\tau = 0$ gives:
$$ \lambda = \alpha = \frac{ \mu^2}{x}  \left({ 1 - 3 x^4}\right)< 0 \Rightarrow x > 4 \mu /3 .$$
Using the fact that $(x^*)^4 = \frac{x^*}{\mu} - 1$, we see that stability of the upper equilibrium is ensured if
\begin{equation} \label{Eqn:xstarcondition} 
 x^* >  \frac{1}{3^{1/4}}  .
 \end{equation}
 
We saw earlier that the maximum of $x - \mu x^4$ occurs at $x_{Max} = \left({ 4 \mu}\right) ^{-1/3}$, and we know that $x^* 
> x_{Max}$.  Because we assume $\mu < \displaystyle{ \frac{3^{3/4}}{4} } $, we see that Equation \ref{Eqn:xstarcondition} is satisfied.  Therefore, when $\tau = 0$, $x^*$ is {\bf stable}.

\begin{figure}
\includegraphics[scale = .4]{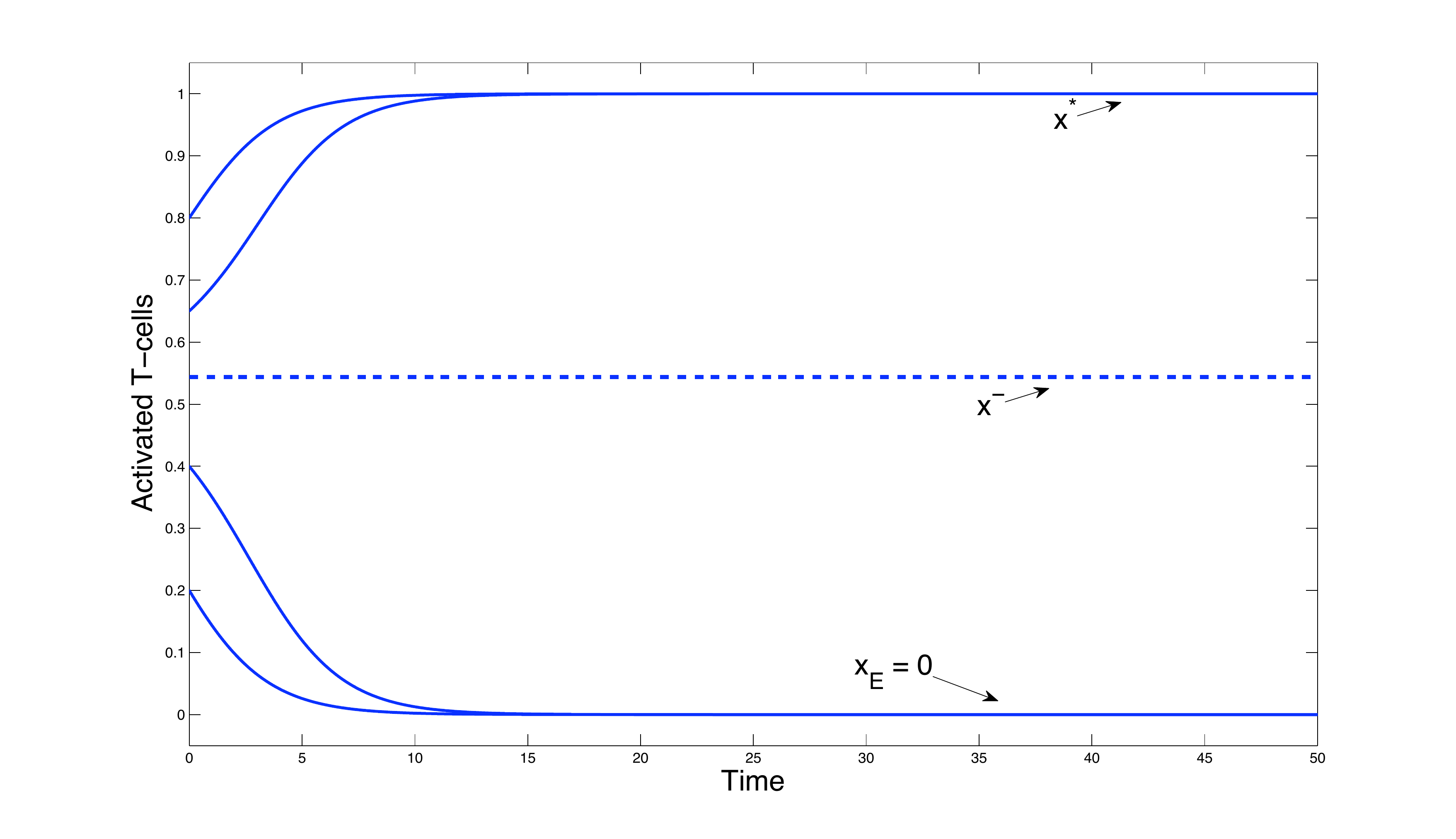}
\caption{Sample trajectories in the case $\tau = 0$ showing two stable equilibria at $x^*$ and $0$, and an unstable equilibrium at $x^-$.  In these simulations $\mu = .5$. } \label{Fig:SampleTrajectories_tauzero}
\end{figure}

Since $\tau = 0$ corresponds to the ODE case,  and $x^-$ lies between the two stable equilibria  $x_E = 0$ and $x_E = x^*$, we know that $x^-$ is {\bf unstable}.  Sample trajectories are shown in Figure \ref{Fig:SampleTrajectories_tauzero}.

\subsection{Stability switching as $\tau$ increases.}

\begin{figure}[h!]
\hskip -.6in \includegraphics[scale = .4]{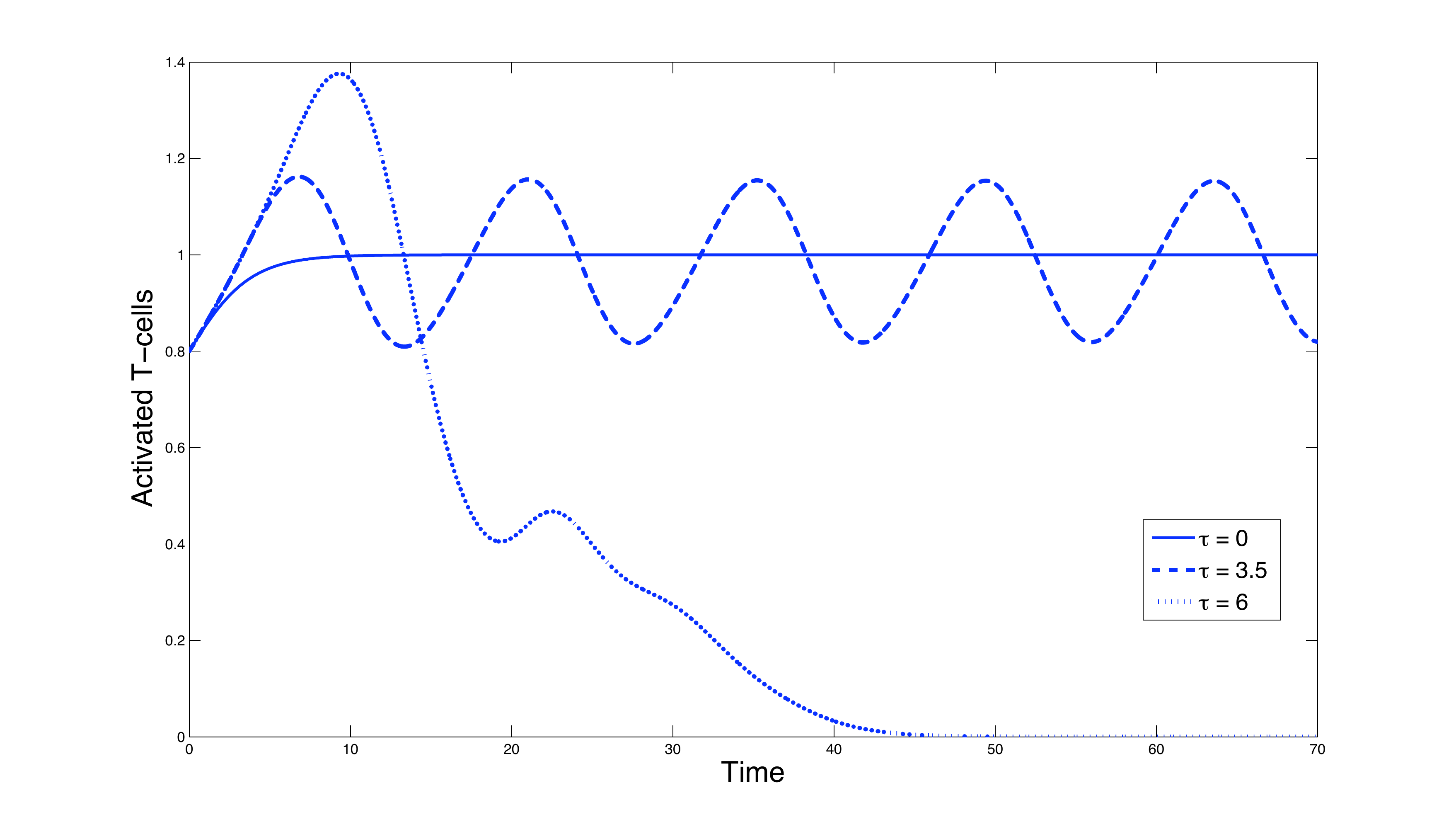}
\caption{As $\tau$ is increased, the equilibrium at $x^* = 1$ loses stability.  First, periodic solutions appear for $\tau$ greater than $\pi$.  Further increases in $\tau$ result in a complete loss of stability, and solutions tend towards the zero equilibrium.}
\label{Fig:SwitchingTimes}
\end{figure}

Changes in stability can occur when
the real part of $\lambda$ switches sign, i.e. when $\lambda$ is purely imaginary.  In this case, the characteristic equation becomes
\begin{equation} \label{Eqn:om}
\lambda = i \omega = \alpha e^{- i \omega \tau} + \beta .
\end{equation}
At the upper equilibrium, $x^*$, $\beta = 0$.  Equating the real and imaginary parts in Equation \ref{Eqn:om} in this case gives:
\begin{eqnarray}
\alpha \cos ( \omega \tau) &=& 0 , \\
- \alpha \sin ( \omega \tau ) &=& \omega .
\end{eqnarray}
The first equation implies that $\omega \tau = (2n + 1) \frac{\pi}{2}$.  Substituting this into the second equation gives:
$$ \omega = \pm \alpha \Rightarrow \tau =  \pm \frac{(2n+1) \pi }{2 \alpha} = \pm \frac{(2n+1) \pi x^*}{2 \mu^2 (1 - 3 {x^*}^4 )} . $$
For a fixed value of $\mu$, this gives a sequence of possible switching times, at which the equilibrium at $x^*$ can change stability.  For example, taking $\mu = \frac{1}{2}$, we see that $x^* = 1$, and possible switching times occur when 
$$ \tau = (2n + 1) \pi, \quad n \in {\mathbb Z} . $$
Figure \ref{Fig:SwitchingTimes} shows numerical solutions for increasing values of $\tau$.  We see that the equilibrium point at 
$x^* = 1$ loses stability when $\tau = \pi$, when a limit cycle appears.  Increasing $\tau$  destroys this limit cycle and, for larger values of $\tau$, solutions that once tended to $x^*$ converge to the zero equilibrium.  

\section{Discussion and Conclusions} \label{Conclusion}

We have presented a simplified model of the interaction between dendritic cells and T-cells at the initiation of a specific immune response.  The model contains one delay that represents the synaptic connection time between dendritic cells and T-cells which is necessary for the development of cytokines that stimulate proliferation.  There is experimental evidence that prolonged contact between dendritic cells and T-cells results in a decrease in viable T-cell production.  The model presented here shows a change in stability of the high T-cell equilibrium as the delay is increased.  This equilibrium goes from being stable when the delay, $\tau$, is small to being unstable for larger values of $\tau$.  These results support the experimental evidence for the ``decreasing potential hypothesis". The model also suggest a role for dendritic cell therapies in boosting the immune defense against certain diseases such as cancer.  However, the model indicates that levels of dendritic cells in the spleen should not be maintained at high levels for too long in order to optimize the expansion of active effector T-cells.

\section*{Acknowledgments} 
We would like to thank the organizers of the 8th AIMS Conference on Dynamical Systems, Differential Equations and Applications that took place in Dresden, Germany, May 25-28, 2010, at which this work was presented.


\end{document}